\documentclass[11pt]{amsart} 
\usepackage{verbatim, latexsym, amssymb, amsmath,color}
\usepackage{epsfig}
\usepackage[colorlinks=true, linkcolor=blue, citecolor=blue]{hyperref}

\def\R{\mathbb R}
\def\RP{\mathbb {RP}}
\def\N{\mathbb N}
\def\Z{\mathbb Z}

\theoremstyle{remark}

\theoremstyle{definition}

\title{Min-max theory and the energy of links}

\author{Ian Agol, Fernando C. Marques and Andr\'e Neves}
\address{%
          University of California, Berkeley\\
970 Evans Hall \#3840\\
Berkeley, CA, 94720-3840}
    \email{%
        ianagol@math.berkeley.edu}  

\address{Instituto de Matem\'atica Pura e Aplicada (IMPA) \\ Estrada Dona Castorina 110 \\ 22460-320 Rio de Janeiro \\ Brazil}
\email{coda@impa.br}
\address{Imperial College \\ Huxley Building \\ 180 Queen's Gate \\ London SW7 2RH \\ United Kingdom}
\email{a.neves@imperial.ac.uk}
\thanks{The first author was supported by DMS-0806027 and DMS-1105738. The second author was partly supported by CNPq-Brazil, FAPERJ, and Math-Amsud. The third author was partly supported by Marie Curie IRG Grant, ERC Start Grant, and Leverhulme Award.}

\begin{document}

\maketitle

\begin{abstract}
Freedman, He, and Wang, conjectured in 1994 that  the M\"{o}bius energy should be minimized, among the class of all nontrivial links in Euclidean space, by the stereographic projection of the standard Hopf link.  We prove this conjecture using the min-max theory of minimal surfaces. 
\end{abstract}

\section{Introduction}

Let $\gamma_i: S^1 \rightarrow \R^3$, $i=1,2$, be a 2-component link, i.e., a pair of rectifiable closed curves in Euclidean three-space  with $\gamma_1(S^1) \cap \gamma_2(S^1) = \emptyset$. The {\it M\"{o}bius cross energy} of the link $(\gamma_1,\gamma_2)$ is defined to be  
$$
E(\gamma_1,\gamma_2) = \int_{S^1 \times S^1} \frac{|\gamma_1'(s)||\gamma_2'(t)|}{|\gamma_1(s)-\gamma_2(t)|^2}\, ds\, dt.
$$ 
The M\"{o}bius energy has the remarkable property of being invariant under conformal transformations of $\R^3$ \cite{freedman-he-wang}. {This allows us to identify freely a link in $\R^3$ with its image in $S^3$ under stereographic projection.}  In the case of knots other energies were considered by O'Hara \cite{o'hara}.

It is not difficult to check that $E(\gamma_1,\gamma_2)\geq 4\pi |{\rm lk}(\gamma_1,\gamma_2)|$, where ${\rm lk}(\gamma_1,\gamma_2)$ denotes the linking number of $(\gamma_1,\gamma_2)$. This is an immediate consequence of the Gauss formula:
$$
{\rm lk}(\gamma_1,\gamma_2) = \frac{1}{4\pi}\int_{S^1 \times S^1} \frac{{\rm det}(\gamma_1'(s),\gamma_2'(t),\gamma_1(s)-\gamma_2(t))}{|\gamma_1(s)-\gamma_2(t)|^3}\, ds \, dt.
$$
By considering pairs of circles which are very far from each other, we see that the cross energy can be made arbitrarily small.  If the linking number of $(\gamma_1,\gamma_2)$ is nonzero, the estimate says that 
$E(\gamma_1,\gamma_2) \geq 4\pi$.  It is natural to search for the optimal configuration in that case.

 It was conjectured by Freedman, He and Wang \cite{freedman-he-wang}, in 1994,  that the M\"{o}bius cross energy should be minimized, among the class of all non-split links in $\R^3$, by the stereographic projection of the standard Hopf link.  
The standard Hopf link $(\hat{\gamma}_1,\hat{\gamma}_2)$ is described by
$$\hat{\gamma}_1(s)=(\cos s, \sin s,0,0) \in S^3\quad\mbox{and}\quad\hat{\gamma}_2(t)=(0,0,\cos t,\sin t) \in S^3,$$ and it is simple to check that  $E(\hat{\gamma}_1,\hat{\gamma}_2)=2\pi^2$. Here we note that the definition of the energy and the conformal invariance property extend 
to any $2$-component link in $\R^n$ \cite{kim-kusner}. A previous result of He implied that the minimal M\"{o}bius cross energy among non-split two component links
is achieved by a link isotopic to the Hopf link \cite[Theorem 4.3]{He02}. It will therefore be a standing assumption throughout the proofs in this paper that $|lk(\gamma_1,\gamma_2)|=1$. 

The goal of this paper is to prove this conjecture:

\subsection{Main Theorem}{\em  Let $\gamma_i: S^1 \rightarrow \R^3$, $i=1,2$, be a 2-component non-split link in $\R^3$.  Then $E(\gamma_1,\gamma_2) \geq 2\pi^2$. 

Moreover, if $E(\gamma_1,\gamma_2)=2\pi^2$ then there exists a conformal map $F:\R^3 \rightarrow {S^3}$ such that $(F\circ \gamma_1,F\circ \gamma_2)=(\hat{\gamma}_1,\hat{\gamma_2})$ describes the standard Hopf link up to orientation and reparameterization.}

\vskip 0.1in
The proof of the Main Theorem follows by applying the min-max theory for the area functional developed  in \cite{marques-neves} to a new five-parameter family of surfaces.  

\subsection{Remark.} The min-max theory was used in \cite{marques-neves} to prove  that the Willmore energy (integral of the square of the mean curvature) of any closed surface in
$\R^3$ with genus $g\geq 1$ is at least $2\pi^2$. This inequality was conjectured by T. J. Willmore in 1965 for the case of tori \cite{willmore}. The Willmore energy of a surface, like the M\"{o}bius cross energy of a link, is  invariant under
conformal transformations. 

\medskip

We now briefly sketch the proof. For any link $(\gamma_1,\gamma_2)$ in $\R^3$, we associate a continuous $5$-parameter family  of surfaces (integral 2-currents with boundary zero, to be more precise)  in $S^3$ such that the area of each surface in the family is bounded above
by $E(\gamma_1,\gamma_2)$. This family is parametrized by a map $\Phi$ defined on $I^5$, and is constructed so that
\begin{itemize}
\item $\Phi(x,0)=\Phi(x,1) = 0$ (trivial surface) for any $x\in I^4$,
\item $\Phi(x,t)$ is an oriented round sphere in $S^3$ for any $x\in \partial I^4$, $t\in [0,1]$,
\item $\{\Phi(x,t)\}_{t\in [0,1]}$ is a homotopically nontrivial sweepout of $S^3$ for any $x\in I^4$,
\item $\sup\{{\rm area}(\Phi(x,t)):(x,t) \in I^5\} \leq E(\gamma_1,\gamma_2)$.
\end{itemize}

This map $\Phi$ has the crucial property that its restriction to $\partial I^4 \times \{1/2\}$ is a homotopically nontrivial map
into  the space of oriented great spheres, which is homeomorphic to $S^3$. Therefore the min-max theory developed  in \cite{marques-neves}  shows the existence of an embedded, smooth, closed  minimal surface $\Sigma \subset S^3$ with genus $g\geq 1$ and such that 
 $${\rm area}\,(\Sigma)\leq \sup\{{\rm area}(\Phi(x,t)):(x,t) \in I^5\}.$$

The fact that $E(\gamma_1,\gamma_2) \geq 2\pi^2$ is a consequence of  the theorem below, proved in \cite{marques-neves}. This theorem rules out the existence of a minimal surface of higher genus with area less than $2\pi^2$:
\subsection*{Theorem}\label{minimal.area.theorem}
\textit{Let $\Sigma\subset S^3$ be an embedded closed minimal surface of genus $g \geq 1$. Then ${\rm area}(\Sigma) \geq 2\pi^2$,
and ${\rm area}(\Sigma)=2\pi^2$ if and only if  $\Sigma$ is the Clifford torus $S^1(\frac{1}{\sqrt{2}})\times S^1(\frac{1}{\sqrt{2}})$ up to  isometries of $S^3$.}

 \medskip

{\bf Acknowledgements:} The authors would like to thank one of the referees for his careful reading and suggestions that helped improving  the exposition of this paper.

\section{Canonical family}

Given a $2$-component link we define a $5$-parameter family of surfaces in $S^3$  with the property that the area of each surface is bounded above by the M\"obius energy of the link. Like in \cite{marques-neves}, we use the language of geometric measure theory, where oriented compact surfaces  become $2$-dimensional integral currents with boundary zero and the area of the surface  is denoted by the mass of the current.  
The point to be careful about is that, due to mass cancellation arising from considering the same set with opposite orientations, the mass is only lower semicontinuous if one considers the weak topology in the space of integral currents. A gentler introduction to the subject can be found in \cite{morgan}.

More precisely, we denote by
\begin{itemize}
\item ${\bf I}_k(S^3)$ the space  of $k$-dimensional integral currents in $\mathbb{R}^4$ with support contained  in $S^3$;
\item ${\mathcal Z}_2(S^3)$ the space  of {integral cycles, i.e.,} integral currents $T \in {\bf I}_2(S^3)$ with  $\partial T=0$.
\end{itemize}
The mass of $T$ is denoted by ${\bf M}(T)$ and it is defined as 
$${\bf M}(T)=\sup\{T(\phi): \phi \in \mathcal{D}^2(\mathbb{R}^4), ||\phi|| \leq 1\}.$$
Here $\mathcal{D}^2(\mathbb{R}^4)$ denotes  the space of smooth $2$-forms in $\mathbb{R}^4$ with compact support, and $||\phi||$ denotes the comass norm of $\phi$.
 
 The {\it flat metric}   on ${\mathcal Z}_2(S^3)$ is defined by 
$$\mathcal{F}(S,T)=\inf\{{\bf M}(P)+{\bf M}(Q): S-T=P+\partial Q, P\in{\bf I}_2(S^3), Q \in {\bf I}_{3}(S^3)\},$$
for $S, T\in {\mathcal Z}_2(S^3)$.  In \cite[Theorem 31.2]{simon} it is shown that, {assuming a uniform bound on the mass of integral cycles,} the topology induced by the flat metric coincides with the weak topology.

Given  any Borel set $A\subset \R^4$, the current $T$ restricted to $A$ is denoted by $T\llcorner A$. Informally, one should think of this as being $T\cap A$.

Given a Lipschitz map $g:S^1\times S^1\rightarrow S^3$ (where $S^1\times S^1$ has a chosen orientation), the current $g_{\#}(S^1\times S^1)\in \mathcal Z_2(S^3)$ is defined by
$$g_{\#}(S^1\times S^1)(\phi)=\int_{S^1\times S^1}g^{*}\phi,\quad\phi\in \mathcal{D}^2(\mathbb{R}^4).$$ We always have
$$
{\bf M}(g_{\#}(S^1\times S^1)) \leq \int_{S^1\times S^1}   |{\rm Jac\, }g|\,ds\,dt.
$$

Let $\gamma_i: S^1 \rightarrow \R^4$, $i=1,2$, be a 2-component link, i.e., a pair of rectifiable curves with $\gamma_1(S^1) \cap \gamma_2(S^1) = \emptyset$. After a reparametrization, we can assume $\gamma_1$ and $\gamma_2$ are Lipschitz and  parametrized {proportional} to  arc length. The M\"{o}bius cross energy $E$ is  invariant under conformal transformations  of $\R^4$ (\cite{kim-kusner}). This means that $E(F\circ \gamma_1,F\circ \gamma_2)=E(\gamma_1,\gamma_2)$ for any conformal map $F:\R^4\rightarrow \R^4$.

The {\it Gauss map} of a link $(\gamma_1,\gamma_2)$ in $\R^4$, denoted by $g=G(\gamma_1,\gamma_2)$, is the Lipschitz map $g:S^1 \times S^1 \rightarrow S^3$ defined by
$$
g(s,t) = \frac{\gamma_1(s)-\gamma_2(t)}{|\gamma_1(s)-\gamma_2(t)|}.
$$
Given an oriented affine hyperplane $P\subset \R^4$, we say that the normal vector $V$ is compatible with the orientation if for every oriented basis $\{e_1,e_2,e_3\}$ of $P$, $\{e_1,e_2,e_3,V\}$ is an oriented basis of $\R^4$. 

We denote an open ball in $\R^4$, centered at $x$ with radius $r$, by $B^4_r(x)$. The boundary of this ball is denoted by $S^3_r(x)$. An intrinsic open ball  of $S^3$,  centered at $p\in S^3$ with radius $r$, is denoted by $B_r(p)$.

\subsection{Lemma}\label{jacobian.estimate}{\it Let $C=g_{\#}(S^1\times S^1)\in  \mathcal{Z}_2(S^3)$. The following properties hold:
\begin{itemize}
\item [(i) ]For almost every $(s,t)\in S^1\times S^1$,
$$
 |{\rm Jac\, }g|(s,t)\leq  \frac{|\gamma_1'(s)||\gamma_2'(t)|}{|\gamma_1(s)-\gamma_2(t)|^2}.
$$ 
If equality holds at $(s,t)$, then 
$$
\langle \gamma_1'(s),\gamma_2'(t) \rangle = \langle \gamma_1'(s),\gamma_1(s)-\gamma_2(t) \rangle=\langle \gamma_2'(t),\gamma_1(s)-\gamma_2(t) \rangle=0.
$$
\item[(ii)]
$${\bf M}(C) \leq \int_{S^1\times S^1} |{\rm Jac\, }g|\,ds\,dt\leq E(\gamma_1,\gamma_2).$$
\item[(iii)] If the link $(\gamma_1,\gamma_2)$ is contained in an oriented affine hyperplane with unit normal vector $p\in S^3$ compatible with the orientation, then
 $$C={\rm lk}(\gamma_1,\gamma_2) \cdot \partial B_{\pi/2}(-p).$$
\end{itemize}
}
\begin{proof}
We have
$$\frac{\partial g}{\partial s}= \frac{1}{|\gamma_1-\gamma_2|}(\gamma_1'-\langle g,\gamma'_1\rangle g)\quad\mbox{and}\quad\frac{\partial g}{\partial t}= -\frac{1}{|\gamma_1-\gamma_2|}(\gamma_2'-\langle g,\gamma'_2\rangle g).$$

Thus
\begin{multline*}\left| \frac{\partial g}{\partial s}\right|^2 \left| \frac{\partial g}{\partial t}\right|^2-\left\langle \frac{\partial g}{\partial s}, \frac{\partial g}{\partial t}\right\rangle^2\leq \left| \frac{\partial g}{\partial s}\right|^2 \left| \frac{\partial g}{\partial t}\right|^2 \\=\frac{|\gamma_1'|^2-\langle g,\gamma'_1\rangle^2}{|\gamma_1-\gamma_2|^2}\frac{|\gamma_2'|^2-\langle g,\gamma'_2\rangle^2}{|\gamma_1-\gamma_2|^2} 
\leq \frac{|\gamma_1'|^2|\gamma_2'|^2}{|\gamma_1-\gamma_2|^4}.
\end{multline*}
This proves the first item.
The second item follows immediately from the first one.

Let $P$ be the oriented hyperplane with normal vector $p\in S^3$ and let $\omega_P$ be its volume form. Let $\omega$, $\omega_{S^3}$, $\omega_{\R^4}$ denote, respectively,   the volume form of $\partial B_{\pi/2}(-p)\subset S^3$ (the exterior unit normal is $p$), $S^3$, and $\R^4$.  We also have
$$\partial B_{\pi/2}(-p)=\{x\in \R^4: \langle x,p\rangle=0\}\cap S^3\subset P.$$

To prove the third item we note that ${\rm supp}\, C\subset \partial B_{\pi/2}(-p)$, and so by  the Constancy Theorem \cite[Theorem 26.27]{simon}  we have $C=k\cdot \partial B_{\pi/2}(-p)$ for some integer $k$. Moreover
\begin{align*}
g^*\omega\left(\frac{\partial}{\partial s},\frac{\partial}{\partial t}\right)& =\omega\left(\frac{\partial g}{\partial s},\frac{\partial g}{\partial t}\right)=\omega_{S^3}\left(\frac{\partial g}{\partial s},\frac{\partial g}{\partial t}, p\right)\\
&=\omega_{\R^4}\left(\frac{\partial g}{\partial s},\frac{\partial g}{\partial t}, p, g\right)= -\omega_{\R^4}\left(\frac{\partial g}{\partial s},\frac{\partial g}{\partial t}, g, p\right)\\
&=-\omega_P\left(\frac{\partial g}{\partial s},\frac{\partial g}{\partial t}, g\right)=
-{\rm det}\left(\frac{\partial g}{\partial s},\frac{\partial g}{\partial t},g\right)\\
& =\frac{{\rm det}(\gamma_1',\gamma_2',\gamma_1-\gamma_2)}{|\gamma_1-\gamma_2|^3},
\end{align*}
and so
$$k=\frac{1}{4\pi} C(\omega)=\frac{1}{4\pi}\int_{S^1\times S^1} g^*\omega={\rm lk}(\gamma_1,\gamma_2).$$

\end{proof}


Given $v\in \R^4$, we define the conformal map $$F_v:\R^4\setminus \{v\} \rightarrow \R^4,\quad F_v(x)=\frac{x-v}{|x-v|^2}.$$ If $v \in B^4$, we have that
$$F_v(S^3_1(0))=S^3_{\frac{1}{1-|v|^2}}(c(v))\quad\mbox{where}\quad c(v)=\frac{v}{1-|v|^2}.$$
\subsection{Lemma}\label{derivative.bound}{\em Given a Lipschitz curve $\gamma:S^1\rightarrow \R^4$ with $|\gamma'(t)| \leq C$, we have
$$\frac{|(F_v\circ\gamma)'(t)|}{|(F_v\circ\gamma)(t)|^2}\leq 3C$$
 for every $v\in \R^4$ and almost every $t\in S^1\setminus \gamma^{-1}(v)$.
}
\begin{proof}
We have
$$(F_v\circ\gamma)'(t)=\frac{1}{|\gamma(t)-v|^2}\left(\gamma'(t)-2\frac{\langle\gamma'(t),\gamma(t)-v\rangle}{|\gamma(t)-v|}\frac{\gamma(t)-v}{|\gamma(t)-v|}\right)$$
for almost every $t\in S^1\setminus \gamma^{-1}(v)$,
and
$$|(F_v\circ\gamma)(t)|^2=\frac{1}{|\gamma(t)-v|^2}.$$
The result follows since $|\gamma'(t)| \leq C$. 

\end{proof}
Given $w\in \R^4$ and $\lambda \in \R$, we set $D_{w,\lambda}(x)=\lambda(x-w)+w$, where $x\in \R^4$.
Finally, given $v\in \overline B^4$ and $z\in (0,1)$, we also define 
 $$b(v,z)=\frac{(2z-1)}{(1-|v|^2+z)(1-z)}$$
 and
 $$a(v,z) =1+(1-|v|^2)b(v,z)= 1+ \frac{(1-|v|^2)(2z-1)}{(1-|v|^2+z)(1-z)}.$$
  For each $v\in B^4$ fixed, $z\to a(v,z)$ is a nondecreasing parametrization of $(0,+\infty)$. 

Suppose now that $\gamma_1(S^1) \cup \gamma_2(S^1) \subset S^3$. 

Given $(v,z)\in B^4\times (0,1)$, we consider $g_{(v,z)}:S^1 \times S^1 \rightarrow S^3$ given by
$$
g_{(v,z)}(s,t) = \frac{(F_v \circ \gamma_1)(s)-(D_{c(v),a(v,z)}\circ F_v \circ \gamma_2)(t)}{|(F_v \circ \gamma_1)(s)-(D_{c(v),a(v,z)}\circ F_v \circ \gamma_2)(t)|},
$$
i.e., $g_{(v,z)}=G(F_v \circ \gamma_1,D_{c(v),a(v,z)}\circ F_v \circ \gamma_2)$.

\subsection{Definition} We define the {\it canonical family} of $(\gamma_1,\gamma_2)$ to be 
 $$C(v,z)={g_{(v,z)}}_\#(S^1\times S^1) \in \mathcal{Z}_2(S^3),$$
 for $(v,z) \in B^4 \times (0,1)$. 
 
Geometrically, $g_{(v,z)}$ corresponds to the Gauss map of the link obtained by applying the conformal transformation $F_v$ to  $(\gamma_1,\gamma_2)$ and then dilating the  curve $F_v\circ\gamma_2$ with respect to the center $c(v)$ by a factor of $a(v,z)$. Both curves $F_v\circ\gamma_1$ and $D_{c(v),a(v,z)}\circ F_v \circ \gamma_2$ are contained in spheres centered at $c(v)$. Notice that  $g_{(v,1/2)}=G(F_v \circ \gamma_1, F_v \circ \gamma_2)$, since $a(v,1/2)=1$.

 Let $\alpha = \inf\{|\gamma_1(s)-\gamma_2(t)|:s,t \in S^1 \}>0$.

\subsection{Lemma}\label{bounds.case3}{\em For all $(v,z)\in B^4\times(0,1)$ we have 
\begin{itemize} 
\item[(i)] 
\begin{eqnarray*}
&&|(F_v \circ \gamma_1)(s)-(D_{c(v),a(v,z)}\circ F_v \circ \gamma_2)(t)|^2\\
&&=a(v,z)|(F_v \circ \gamma_1)(s)-(F_v \circ \gamma_2)(t)|^2+b(v,z)^2;
\end{eqnarray*}
\item[(ii)] 
$$|(F_v \circ \gamma_1)(s)-(F_v \circ \gamma_2)(t)|^2 \geq \alpha^2 |(F_v \circ \gamma_1)(s)|^2|(F_v \circ \gamma_2)(t)|^2$$
for all $s,t\in S^1$.
\end{itemize}
}
\begin{proof}
 If $x,y \in S^3_r(p)$, then
$$
|x - (p+a(y-p))|^2 = a\, |x-y|^2+(1-a)^2r^2.
$$
The identity of item (i) follows because $(F_v \circ \gamma_1)(s), (F_v \circ \gamma_2)(t) \in S_{r}(p)$ where  $r=(1-|v|^2)^{-1}$ and $p=(1-|v|^2)^{-1}v=c(v)$.

The second item follows from
\begin{eqnarray*}
|(F_v \circ \gamma_1)(s)-(F_v \circ \gamma_2)(t)|^2 &=& \left|\frac{\gamma_1(s)-v}{|\gamma_1(s)-v|^2}-\frac{\gamma_2(t)-v}{|\gamma_2(t)-v|^2}\right|^2\\
&=& \frac{1}{|\gamma_1(s)-v|^2}-2\frac{\langle \gamma_1(s)-v,\gamma_2(t)-v\rangle}{|\gamma_1(s)-v|^2|\gamma_2(t)-v|^2}\\
& &+\frac{1}{|\gamma_2(t)-v|^2}\\
&=& \frac{|\gamma_1(s)-\gamma_2(t)|^2}{|\gamma_1(s)-v|^2|\gamma_2(t)-v|^2}\\
&\geq&\frac{\alpha^2}{|\gamma_1(s)-v|^2|\gamma_2(t)-v|^2}\\
& = & \alpha^2|(F_v \circ \gamma_1)(s)|^2|(F_v \circ \gamma_2)(t)|^2.
\end{eqnarray*}

\end{proof}
 
 The following  lemma establishes  an important property of the canonical family.
 
 \subsection{Lemma}\label{canonical.area.bound}{\em  For every $(v,z)\in B^4\times (0,1)$, we have
\begin{multline*}
\int_{S^1\times S^1}|{\rm Jac}\,g_{(v,z)}|ds\,dt \\
\leq   \int_{S^1\times S^1}\frac{a(v,z)|(F_{v}\circ \gamma_1)'(s)||(F_{v}\circ \gamma_2)'(t)|}{a(v,z)|F_{v}\circ \gamma_1(s)-F_{v}\circ \gamma_2(t)|^2+b(v,z)^2}ds\,dt\\
\leq E(\gamma_1,\gamma_2).
\end{multline*}
In particular, ${\bf M}(C(v,z))\leq   E(\gamma_1,\gamma_2).$}
\begin{proof}
From Lemma \ref{bounds.case3} (i)  we have 
	$$|F_{v}\circ \gamma_1-D_{c(v), a(v,z)}\circ F_{v}\circ \gamma_2|^2= a(v,z)|F_{v}\circ \gamma_1-F_{v}\circ \gamma_2|^2+b(v,z)^2.$$
	Moreover 
	$$|(D_{c(v), a(v,z)}\circ F_{v}\circ \gamma_2)'(t)|=a(v,z)|(F_{v}\circ \gamma_2)'(t)|.$$
	Thus, combining with  Lemma \ref{jacobian.estimate} (i) we obtain for almost all $(s,t)\in S^1\times S^1$

	\begin{multline*}
	|{\rm Jac}\,g_{(v,z)}|(s,t)\leq \frac{|(F_{v}\circ \gamma_1)'(s)| |(D_{c(v), a(v,z)}\circ F_{v}\circ \gamma_2)'(t)|}{|F_{v}\circ \gamma_1(s)-D_{c(v), a(v,z)}\circ F_{v}\circ \gamma_2(t)|^2}\\
	{=}\frac{a(v,z)|(F_{v}\circ \gamma_1)'(s)||(F_{v}\circ \gamma_2)'(t)|}{a(v,z)|F_{v}\circ \gamma_1(s)-F_{v}\circ \gamma_2(t)|^2+b(v,z)^2}.
	\end{multline*}
	Integrating the above {inequality} gives the desired result.
	\end{proof}

We also get  the following uniform control on the maps $g_{(v,z)}$.

\subsection{Lemma}\label{energy.bound}{\em   There exists a constant $C>0$ {(depending on $(\gamma_1,\gamma_2)$)} such that for all $(v,z)\in B^4\times (0,1)$ and almost all $(s,t)\in S^1\times S^1$  we have
$$|{\rm Jac}\,g_{(v,z)}|(s,t)\leq C.$$
}
\begin{proof}

	Combining Lemma \ref{bounds.case3} (i) and (ii) we have 
	$$|F_{v}\circ \gamma_1-D_{c(v), a(v,z)}\circ F_{v}\circ \gamma_2|^2\geq \alpha^2 a(v,z)|F_{v}\circ \gamma_1|^2|F_{v}\circ \gamma_2|^2.$$
	Moreover 
	$$|(D_{c(v), a(v,z)}\circ F_{v}\circ \gamma_2)'(t)|=a(v,z)|(F_{v}\circ \gamma_2)'(t)|.$$
	Combining both inequalities we have from Lemma \ref{derivative.bound} the existence of some constant $C>0$ so that 
	\begin{multline*}
	\frac{|(F_{v}\circ \gamma_1)'(s)| |(D_{c(v), a(v,z)}\circ F_{v}\circ \gamma_2)'(t)|}{|F_{v}\circ \gamma_1(s)-D_{c(v), a(v,z)}\circ F_{v}\circ \gamma_2(t)|^2}\\
	\leq \frac{1}{\alpha^2}\frac{|(F_{v}\circ \gamma_1)'(s)|}{|F_{v}\circ \gamma_1(s)|^2}\frac{|(F_{v}\circ \gamma_2)'(t)|}{|F_{v}\circ \gamma_2(t)|^2}
	\leq {C} \,  a.e.
	\end{multline*}
	Thus, from Lemma \ref{jacobian.estimate} (i), we obtain that for almost all $(s,t)\in S^1\times S^1$
	 $$|{\rm Jac}\,g_{(v,z)}|(s,t)\leq {C}.$$
	\end{proof}


\section{Continuity of canonical family}\label{continuity}
In this section we show that the continuous map  given by the canonical family $$C:B^4\times (0,1)\rightarrow \mathcal{Z}_2(S^3)$$  can be continuously extended to $\overline B^4\times[0,1]$. In order to do that, we first  extend the definition of $g_{(v,z)}$ when $(v,z)\in \overline B^4\times[0,1]$.

 If $v\in S^3$ and {$0<z<1$}, consider the hyperplane
$$P_{(v,z)}=\{x\in\R^4\,:\,\langle x,v\rangle = -1/2-b(z)\},\quad\mbox{where}\quad b(z)=\frac{2z-1}{z(1-z)}.$$
Define, for every $0<z<1$, the smooth map $$L_{(v,z)}:S^3\setminus\{v\}\rightarrow P_{(v,z)},\quad L_{(v,z)}(x)=F_v(x)-b(z)v.$$
 Note that indeed
$$
\langle L_{(v,z)}(x),v\rangle=\left\langle \frac{x-v}{|x-v|^2}, v \right\rangle -b(z) = \frac{\langle x,v\rangle -1}{2-2\langle x,v\rangle} -b(z)=-1/2-b(z).
$$
{We shall see in the next proposition that if $\{v_i\}_{i\in\N}$ in $B^4$ tends to $v\in S^3$ as $i$ tends to infinity then, for all $x\in  S^3\setminus\{v\}$, $D_{c(v_i),a(v_i,z)}\circ F_{v_i}(x)$ tends to $L_{(v,z)}(x)$.}

Set 
  $$S_v= (\gamma_1^{-1}(v) \times S^1) \cup (S^1 \times  \gamma_2^{-1}(v)).$$
  If $v\in S^3$ and $0<z<1$ we  define
$$g_{(v,z)}:(S^1\times S^1) \setminus S_v \rightarrow S^3, \quad g_{(v,z)}=G(F_v\circ\gamma_1, L_{(v,z)}\circ \gamma_2).$$
 For $(v,z)\in \overline B^4\times\{0,1\}$ we define, for all $x\in(S^1\times S^1) \setminus S_v$, {$g_{(v,1)}(x)=v$ and $g_{(v,0)}(x)=-v$.}

\subsection{Proposition}\label{equicontinuous.extension}{\em The canonical family  can be  extended to a continuous map in the flat metric $$C:\overline B^4\times [0,1]\rightarrow \mathcal{Z}_2(S^3)$$ with the following properties:
\begin{itemize}
\item[(i)] $C(v,0)=C(v,1)=0$ for all $v\in \overline B^4$.
\item[(ii)]
If $v\in S^3\setminus(\gamma_1(S^1)\cup \gamma_2(S^1))$ and $0<z<1$ then 
$$C(v,z)={g_{(v,z)}}_{\#}(S^1\times S^1).$$ 
\item[(iii)] If $(v_i,z_i)\in B^4\times (0,1)$ tends to $(v,z)\in \overline B^4\times[0,1]$, then
$$\lim_{i\to\infty} \int_{S^1\times S^1}  |{\rm Jac}\,g_{(v_i,z_i)}| \, ds\, dt= \int_{(S^1\times S^1)\setminus S_v}  |{\rm Jac}\,g_{(v,z)}| \, ds\, dt.$$
\item[(iv)]
$$\lim_{r\to 0}\sup \left\{{\bf M}(C(v,z)\llcorner B_r(p)):\,p\in S^3, (v,z)\in B^4\times (0,1) \right\} =0.$$
\item[(v)] For every $v\in S^3$, $C(v,1/2)=-{\rm lk}(\gamma_1,\gamma_2)\cdot\partial B_{\pi/2}(v)$.
\end{itemize}
}
\subsection{Remark:}   While the mass is only lower semicontinuous for the weak topology, Proposition \ref{equicontinuous.extension} (iii) implies the $L^1$ norm of the Jacobian is continuous. This is the reason we look at this quantity  instead of the mass.

Proposition \ref{equicontinuous.extension} (iv) is a technical condition needed  to apply the work in \cite{marques-neves} and it says that no $C(v,z)$ has large  mass concentrated in a small ball.
\vskip 0.1in
 \begin{proof}
 Since  $\gamma_1$ and $\gamma_2$ are both Lipschitz and parametrized proportionally to arc length, $\mathcal H^1(\gamma_j^{-1}(v))=0$  (\cite[Corollary 2.7.5]{burago-burago-ivanov}) for $j=1,2$ and any $v\in \overline{B}^4$. Thus $S_v$ is a closed set with $\mathcal{H}^2(S_v)=0$ (where $\mathcal{H}^i$ denotes $i$-dimensional Hausdorff measure). 
 
Consider a sequence $(v_i,z_i)\in B^4\times(0,1)$ tending to $(v,z)\in \overline B^4\times [0,1]$.

 We first argue that, when $0<z<1$,
 $$\lim_{i\to\infty}g_{(v_i,z_i)}(s,t)=g_{(v,z)}(s,t)\quad\mbox{for all}\quad  (s,t)\in (S^1\times S^1)\setminus S_v.$$
 If $v\in B^4$ this follows directly  from the definition of the maps $g_{(v,z)}$ (note that in this case $S_v=\emptyset$). If $v\in S^3$ then, for every $x\in S^3\setminus\{v\}$, we have
\begin{align*}
\lim_{i\to\infty}D_{c(v_i),a(v_i,z)}\circ F_{v_i}(x) & 
=\lim_{i\to\infty}a(v_i,z)F_{v_i}(x)-v_i\frac{a(v_i,z)-1}{1-|v_i|^2}\\
& =F_v(x)-b(z)v=L_{(v,z)}(x).
\end{align*} 
 Thus $g_{(v_i,z_i)}$ converges pointwise to the Lipschitz map $g_{(v,z)}$ on  $ (S^1\times S^1)\setminus S_v.$ Similarly we can prove that the first derivative $Dg_{(v_i,z_i)}$ converges pointwise to  $Dg_{(v,z)}$ almost everywhere on  $ (S^1\times S^1)\setminus S_v.$

From Lemma \ref{energy.bound} we have that $ |{\rm Jac}\,g_{(v_i,z_i)}|$ is uniformly bounded and thus, combining the previous  pointwise convergences  with Lebesgue's Dominated Convergence Theorem, we have that for all $\phi\in \mathcal{D}^2(\R^4)$
\begin{equation}\label{unique.limit1}
\lim_{i\to\infty} C(v_i,z_i)(\phi)=\lim_{i\to\infty}\int_{S^1\times S^1}g_{(v_i,z_i)}^*(\phi)= \int_{(S^1\times S^1)\setminus S_v}g_{(v,z)}^*(\phi),
\end{equation}
and
\begin{equation}\label{unique.limit2}
\lim_{i\to\infty}\int_{S^1\times S^1} |{\rm Jac}\,g_{(v_i,z_i)}| \, ds\,dt= \int_{(S^1\times S^1)\setminus S_v}|{\rm Jac}\,g_{(v,z)}| \, ds\,dt.
\end{equation}

We now extend the two identities above to the case when $z=0$ or $z=1$. It suffices to argue that
\begin{equation*}
\lim_{i\to\infty}\int_{S^1\times S^1} |{\rm Jac}\,g_{(v_i,z_i)}|ds\, dt=0.
\end{equation*}

We have  from Lemma \ref{bounds.case3} (i) and (ii) that
	\begin{multline*}
	|F_{v_i}\circ \gamma_1-D_{c(v_i), a(v_i,z_i)}\circ F_{v_i}\circ \gamma_2|^2\\
	\geq  {\alpha^2} {a(v_i,z_i)}|F_{v_i}\circ \gamma_1|^2|F_{v_i}\circ \gamma_2|^2+b(v_i,z_i)^2.
	\end{multline*}
	Moreover, 
	$$|(D_{c(v_i), a(v_i,z_i)}\circ F_{v_i}\circ \gamma_2)'|=a(v_i,z_i)|(F_{v_i}\circ \gamma_2)'|$$
and thus, from Lemma \ref{jacobian.estimate} (i) and Lemma \ref{derivative.bound}, we have  for some constant $C>0$ that
\begin{multline*}
	|{\rm Jac}\, g_{(v_i,z_i)}|\leq \frac{|(F_{v_i}\circ \gamma_1)'| |(D_{c(v_i), a(v_i,z_i)}\circ F_{v_i}\circ \gamma_2)'|}{|F_{v_i}\circ \gamma_1-D_{c(v_i), a(v_i,z_i)}\circ F_{v_i}\circ \gamma_2|^2}\\
	\leq \frac{a(v_i,z_i)|(F_{v_i}\circ \gamma_1)'||(F_{v_i}\circ \gamma_2)'|}{ {\alpha^2} {a(v_i,z_i)}|F_{v_i}\circ \gamma_1|^2|F_{v_i}\circ \gamma_2|^2+b(v_i,z_i)^2}\\
	\leq \frac{Ca(v_i,z_i)|F_{v_i}\circ \gamma_1|^2|F_{v_i}\circ \gamma_2|^2}{{\alpha^2} {a(v_i,z_i)}|F_{v_i}\circ \gamma_1|^2|F_{v_i}\circ \gamma_2|^2+b(v_i,z_i)^2}.
	\end{multline*}
Given $(s,t)\in (S^1\times S^1)\setminus S_v$, we have that $|F_{v_i}\circ \gamma_1|(s)$ and $|F_{v_i}\circ \gamma_2|(t)$ are uniformly bounded. 

We claim that 
$a(v_i,z_i)/b(v_i,z_i)^2\to 0$. This implies that $|{\rm Jac}\,g_{(v_i,z_i)}|$ converges pointwise to zero on $(S^1\times S^1)\setminus S_v$ and so we conclude from the  Lebesgue's Dominated Convergence Theorem
$$ \lim_{i\to\infty}\int_{S^1\times S^1} |{\rm Jac}\,g_{(v_i,z_i)}|\, ds\,dt=0.
$$
If $z_i\to 0$ and $v\in B^4$, then $a(v_i,z_i)\to 0$, $b(v_i,z_i)\to (|v|^2-1)^{-1}$ and so indeed $a(v_i,z_i)/b(v_i,z_i)^2$ tends to zero.
If {$v\in S^3$ and either $z_i\to 1$ or $z_i\to 0$}, then 
  $|b(v_i,z_i)|\to +\infty$ and we have 
$$
\lim_{i\to\infty} \frac{a(v_i,z_i)}{b^2(v_i,z_i)}=\lim_{i\to\infty} \frac{1}{b^2(v_i,z_i)}+\frac{1-|v_i|^2}{b(v_i,z_i)}=0.
$$

We now have all the necessary estimates to prove Proposition \ref{equicontinuous.extension}.

From  Federer-Fleming Compactness Theorem \cite[Theorem 27.3]{simon},  we know that $C(v_i,z_i)$ has a convergent subsequence. It follows from   \eqref{unique.limit1} that the limit of this convergent subsequence is uniquely determined. The fact that the sequence $(v_i,z_i)$ was arbitrary implies that the map $C$ admits  a continuous extension to $\overline B^4\times [0,1].$  In case $v\in S^3\setminus(\gamma_1(S^1)\cup \gamma_2(S^1))$, the set $S_v$ is empty and so Proposition \ref{equicontinuous.extension} (ii) follows from \eqref{unique.limit1}.

Proposition \ref{equicontinuous.extension} (i) and (iii) follow from  \eqref{unique.limit2}.

We now prove property (iv). It suffices to show that for every $\delta>0$ and $p\in S^3$, we can find $r=r(p,\delta)$ so that 
$$\int_{g^{-1}_{(v,z)}(B_r(p))} |{\rm Jac}\,g_{(v,z)}| \, ds\, dt\leq \delta\quad\mbox{for all}\quad (v,z)\in B^4\times (0,1).$$
This is because, via a standard finite covering argument, we can then find $r$ independent of $p$.

Suppose this statement is false. This means that there exist $p\in S^3$, $\delta>0$, and a sequence $(v_i,z_i) \in B^4 \times (0,1)$ such that 
$$\int_{g^{-1}_{(v_i,z_i)}(B_{1/i}(p))}|{\rm Jac}\,g_{(v_i,z_i)}|\, ds\, dt\geq \delta$$
for every $i \in \N$. By passing to a subsequence we can assume $(v_i,z_i)$ converges to $(v,z) \in \overline{B}^4 \times [0,1]$. 

Because $|{\rm Jac}\,g_{(v_i,z_i)}|$ is uniformly  bounded and converges almost everywhere to $|{\rm Jac}\,g_{(v,z)}|$, we have
$$\int_{g_{(v,z)}^{-1}(p)\cap\left((S^1\times S^1)\setminus S_v \right)}|{\rm Jac}\,g_{(v,z)}| ds\, dt\geq \delta .$$
The area formula \cite[Identity 8.4]{simon} shows the inequality above is impossible. This is a contradiction, thus property (iv) holds.

To show property (v), note that given $v\in S^3\setminus(\gamma_1(S^1)\cup \gamma_2(S^1))$, $$(F_v\circ \gamma_1, F_v\circ \gamma_2)$$ is a link in the affine hyperplane $P_{(v,1/2)}$, where $F_v$ sends the exterior unit normal of $S^3$ into $v$. Thus from Lemma \ref{jacobian.estimate} (iii)  we have
\begin{multline*}
C(v,1/2)=G(F_v\circ\gamma_1, F_v\circ \gamma_2)_{\#}(S^1\times S^1)\\
={\rm lk}(\gamma_1,\gamma_2)\cdot\partial B_{\pi/2}(-v)=-{\rm lk}(\gamma_1,\gamma_2)\cdot\partial B_{\pi/2}(v).
\end{multline*}
The continuity of $C$ allows us to extend the above identity to every $v\in S^3$.
  \end{proof}

\section{Extension of family}
 
The main arguments in \cite{marques-neves} relied in a crucial way on the fact that the restriction of the  $5$-parameter family, defined on $\overline B^4\times[0,1]$,  to $S^3\times[0,1]$   was a non-trivial map into the set of all  oriented round spheres of $S^3$. The canonical family defined in the previous section does not have this property.

In this section we construct a continuous extension of $C$ to $\overline{B}^4_2(0) \times [0,1]$ such that its restriction  to $\partial B^4_2(0) \times [0,1]$ is a non-trivial map into the space of oriented geodesic spheres of $S^3$. When doing this, it is important for the proof of the main theorem that the  mass of each member of this new family is still bounded above by $E(\gamma_1,\gamma_2)$.

To achieve this property we use the fact that, for each $p\in S^3$, the support of the current $C(p,z)$ lies in the northern hemisphere $\{x:\langle x,p\rangle>0\}$ if $z>1/2$ and in the southern hemisphere  $\{x:\langle x,p\rangle<0\}$ if $z<1/2$. This allows us to continuously deform each $C(p,z)$ into a geodesic sphere without increasing the mass.

 \subsection{Proposition} {\it There exists a constant $c>0$ such that for every $p \in S^3$ we have
 \begin{itemize}
 \item ${\rm supp}(C(p,z)) \subset\overline{B}_{\pi/2}(p) \setminus B_{r(z)}(p)$ if $z\in [1/2,1]$,
 \item ${\rm supp}(C(p,z)) \subset \overline{B}_{\pi/2}(-p) \setminus B_{\pi-r(z)}(-p)$ if $z\in [0,1/2]$,
\end{itemize}
where $$r(z)=\cos^{-1} \left(\frac{b(z)}{\sqrt{|b(z)|^2 + c^2}}\right)\in [0,\pi]\quad\mbox{and}\quad b(z)=\frac{2z-1}{z(1-z)}.$$}

\medskip

Notice that $r(0)=\pi$, $r(1/2)=\pi/2$,  and  $r(1)=0$.

\begin{proof}
If $z=0$ or $z=1$, we have $C(p,z)=0$ by Theorem \ref{equicontinuous.extension} (i) and the proposition follows immediately.

Suppose $z\in (0,1)$ and $p\in S^3 \setminus (\gamma_1(S^1) \cup \gamma_2(S^1))$. 

From Theorem \ref{equicontinuous.extension} (ii) we have that ${\rm supp}(C(p,z)) \subset g_{(p,z)}(S^1\times S^1)$. But
\begin{eqnarray*}
\langle g_{(p,z)}, p \rangle &=& \left\langle \frac{F_{p}\circ \gamma_1-L_{(p,z)}\circ \gamma_2}{|F_{p}\circ \gamma_1-L_{(p,z)}\circ \gamma_2|}, p  \right\rangle\\
&=&\left\langle \frac{F_{p}\circ \gamma_1-F_{p}\circ \gamma_2+b(z)p}{|F_{p}\circ \gamma_1-L_{(p,z)}\circ \gamma_2|},p \right \rangle\\
&=& \frac{b(z)}{|F_{p}\circ \gamma_1-L_{(p,z)}\circ \gamma_2|}.\\
\end{eqnarray*}

This already implies that $\langle g_{(p,z)}, p \rangle \geq 0$ if $z\in [1/2,1]$, and that $\langle g_{(p,z)}, p \rangle \leq 0$ if $z\in [0,1/2]$. 

Using the fact that for some constant $c>0$ we have
$$|F_p\circ\gamma_1(s)- F_p\circ\gamma_2(t)|^2\geq c^2=\alpha^2/16$$
from the argument in Lemma \ref{bounds.case3}, 
for all $p\in S^3$ and  $(s,t)\in S^1\times S^1$,  it  follows that
\begin{eqnarray*}
|\langle g_{(p,z)}(s,t), p \rangle| &\leq&  \frac{|b(z)|}{\sqrt{c^2+ b(z)^2}}.\\
\end{eqnarray*}
This proves the proposition for $z\in [0,1]$ and $p\in S^3 \setminus (\gamma_1(S^1) \cup \gamma_2(S^1))$. 

Since $S^3 \setminus (\gamma_1(S^1) \cup \gamma_2(S^1))$ is everywhere dense in $S^3$, and $C:S^3 \times [0,1] \rightarrow \mathcal{Z}_2(S^3)$ is continuous in the flat topology, the proposition also holds for $p \in \gamma_1(S^1) \cup \gamma_2(S^1)$. 
\end{proof}

 Let $p\in S^3$, $\lambda \in [0,\pi/2]$. For $t\in [0,1]$, we define the retraction map $$R_{(p,\lambda,t)}:\overline{B}_{\pi/2}(p)\setminus B_\lambda(p) \rightarrow \overline{B}_{\pi/2}(p)\setminus B_\lambda(p)$$ by
 $$
 R_{(p,\lambda,t)}(x) =  {\rm {exp}}_p\left( \left((1-t)+t\frac{\lambda}{d(p,x)}\right){\rm exp}_p^{-1}(x)\right).
 $$
 
 Notice that $R_{(p,\pi/2,t)}:\partial B_{\pi/2}(p) \rightarrow \partial B_{\pi/2}(p)$ is the identity map for every $p\in S^3$ and $t\in [0,1]$.

\subsection{Proposition}\label{retraction}{\it We have
\begin{itemize}
 \item $R_{(p,\lambda,0)}(x)=x$ for all $x \in \overline{B}_{\pi/2}(p)\setminus B_\lambda(p)$;

\item $R_{(p,\lambda,1)}\left(\overline{B}_{\pi/2}(p)\setminus B_\lambda(p) \right) \subset \partial B_\lambda(p)$;

 \item $|(DR_{(p,\lambda,t)})_x(v)| \leq |v|$ for all $x \in \overline{B}_{\pi/2}(p)\setminus B_\lambda(p)$, $v\in T_xS^3$.
 \end{itemize}
 }
 
 \begin{proof}
It follows from the definition that $$R_{(p,\lambda,0)}(x)={\rm exp}_p\left( {\rm exp}_p^{-1}(x)\right)=x,$$ and that
$$
R_{(p,\lambda,1)}(x) =  {\rm exp}_p\left( \frac{\lambda}{d(p,x)}{\rm exp}_p^{-1}(x)\right).
$$
Since $|{\rm exp}_p^{-1}(x)|=d(p,x)$, the equalities above prove the first two items of the proposition.

Note that ${\rm exp}_p^{*}(\overline{g})_y = dr^2 + \sin^2(r) d\omega^2$, where $\overline{g}$ denotes the standard metric on $S^3$, $d\omega^2$ denotes the standard metric on $S^2$, and $r(y)=|y|=d(p,{\rm exp}_p(y))$. 

The map 
$f={\rm exp}_p^{-1}\circ R_{(p,\lambda,t)}\circ {\rm exp}_p$ is given by
$$
f(y) = \left((1-t)+t\frac{\lambda}{|y|}\right) y.
$$
Hence $Df_y(\partial_r) = (1-t) \partial_r$, where $\partial_r=y/|y|$, and $$Df_y(v)=\left((1-t)+t\frac{\lambda}{|y|}\right) v$$ if $\langle v, \partial_r\rangle =0$.

Thus
$$
{\rm exp}_p^{*}(\overline{g})_{f(y)}(Df_y(\partial_r),Df_y(\partial_r)) = (1-t)^2 \leq 1={\rm exp}_p^{*}(\overline{g})_{y}(\partial_r,\partial_r)
$$
and
$${\rm exp}_p^{*}(\overline{g})_{f(y)}(Df_y(\partial_r),Df_y(v))=0={\rm exp}_p^{*}(\overline{g})_{y}(\partial_r,v),\quad\mbox{if}\quad \langle v, \partial_r\rangle =0.
$$

 The third item of this proposition follows since, if $\langle v, \partial_r\rangle =0$ and $\lambda \leq |y| \leq \pi/2$, we have
 \begin{eqnarray*}
 {\rm exp}_p^{*}(\overline{g})_{f(y)}(Df_y(v),Df_y(v)) &=& \frac{\sin^2(|f(y)|)}{|f(y)|^2}\left((1-t)+t\frac{\lambda}{|y|}\right) ^2 |v|_{Eucl}^2\\
 &=&  \frac{\sin^2(|f(y)|)}{|y|^2} |v|_{Eucl}^2\\
 &\leq& \frac{\sin^2(|y|)}{|y|^2} |v|_{Eucl}^2\\
 &=& {\rm exp}_p^{*}(\overline{g})_{y}(v,v).
 \end{eqnarray*}
 We have used the fact that $r \rightarrow \sin r$ is nondecreasing in $[0,\pi/2]$, and that $|f(y)|\leq |y|$ whenever
 $\lambda \leq |y|$.
 \end{proof}

For $p\in S^3$, $z\in [0,1]$, $s\in [0,1]$, we set $C(p,z,s) \in \mathcal{Z}_2(S^3)$ by
\begin{eqnarray*}
C(p,z,s) = 
\left\{
\begin{array}{ll}
{R_{(p,r(z),s)}}_\# C(p,z)  & {\rm  \ if \ } z\in [1/2,1] \\
{R_{(-p,\pi-r(z),s)}}_\# C(p,z)  & {\rm \ if \ } z\in [0,1/2]. \\
\end{array}
\right.
\end{eqnarray*}

\subsection{Proposition}\label{mapa.K}{\it Assume ${\rm lk}(\gamma_1,\gamma_2)=-1$. The map $$C:S^3 \times [0,1] \times [0,1] \rightarrow \mathcal{Z}_2(S^3)$$ satisfies the following properties:
\begin{itemize}
\item[(i)] $C$ is continuous in the flat topology;
\item[(ii)] $C(p,z,0)=C(p,z)$ for all $p\in S^3$, $z\in [0,1]$;
\item[(iii)] $C(p,z,1)=\partial B_{r(z)}(p)$ for all $p\in S^3$, $z\in [0,1]$;
\item[(iv)] ${\bf M}(C(p,z,s)) \leq {\bf M}(C(p,z))$ for all $p\in S^3$, $z\in [0,1]$, $s\in [0,1]$;
\item[(v)] $$
\lim_{r\to 0}\sup \left\{{\bf M}(C(p,z,s)\llcorner B_r(q))\,:\,p, q\in S^3\mbox{ and } z,s\in [0,1]\right\}=0. 
$$
\end{itemize}

}

\begin{proof}
The  family  $(p,\lambda,t)\mapsto R_{(p,\lambda,t)}$ is continuous in the space of $C^1$ maps and so, since $(p,z) \mapsto C(p,z)$
is already continuous in the flat topology, item (i) follows. 

Item (ii) follows immediately from the fact that $R_{(p,\lambda,0)}=id$.

Consider $\hat{C}:S^3 \times [0,1] \rightarrow \mathcal{Z}_2(S^3)$ defined by $\hat{C}(p,z)=C(p,z,1)$. Then $\hat{C}$ is continuous in the flat topology. Since $${\rm supp}(\hat{C}(p,z))\subset R_{(p,r(z),1)}({\rm supp}(C(p,z)))\subset \partial B_{r(z)}(p),$$
and $\partial \hat{C}(p,z)=0$, the Constancy Theorem \cite[Theorem 26.27]{simon} implies that there exists $k(p,z)\in \mathbb{Z}$ such that
$$
\hat{C}(p,z) = k(p,z)\cdot \partial B_{r(z)}(p) \in \mathcal{Z}_2(S^3).
$$

Because of the continuity of $\hat{C}$  in the flat topology, there must exist $k\in \mathbb{Z}$ such that $k(p,z)=k$ for every $p\in S^3$ and $z\in [0,1]$. From Theorem \ref{equicontinuous.extension} (v)
\begin{eqnarray*}
\hat{C}(p,1/2)&=&C(p,1/2,1)\\
&=&{R_{(p,r(1/2),1)}}_\# (C(p,1/2)) \\
&=&{R_{(p,\pi/2,1)}}_\# (\partial B_{\pi/2}(p)) \\
&=& id_\# (\partial B_{\pi/2}(p))  =  \partial B_{\pi/2}(p).
\end{eqnarray*}
Hence $k=1$, and this proves item (iii).

If $f$ is a Lipschitz map with $|df| \leq 1$ and $C$ is an integral current, then ${\bf M}(f_\#C) \leq {\bf M}(C)$ \cite[Lemma 26.25]{simon}. Hence item (iv) follows from the definition of $C(p,z,s)$ and Proposition \ref{retraction}.

 Given $p\in S^3$, $0<z<1$, and $0\leq s\leq 1$, set $g_{(p,z,s)}=R_{(p,r(z),s)}\circ g_{(p,z)}.$  
 From Proposition \ref{retraction} we have that $|{\rm Jac}\,g_{(p,z,s)}|\leq |{\rm Jac}\,g_{(p,z)}|$.  We also have that for every sequence $(p_i,z_i,s_i)$ tending to $(p,z,s)$, $|{\rm Jac}\,g_{(p_i,z_i,s_i)}|$ converges pointwise to $|{\rm Jac}\,g_{(p,z,s)}|$ in $S^1\times S^1\setminus S_p$. We can then use Lebesgue's Dominated Convergence Theorem and the {area} formula, exactly like in  the proof of Proposition \ref{equicontinuous.extension} (iv), to conclude item (v).
\end{proof}

We now define the extension $\tilde{C}:\overline{B}^4_2(0)\times [0,1] \rightarrow \mathcal{Z}_2(S^3)$ by
\begin{eqnarray*}
\tilde{C}(v,t) = 
\left\{
\begin{array}{ll}
C(v,t) & {\rm \ if \ } v \in \overline{B}^4_1(0) \\
C\left(\frac{v}{|v|},t,|v|-1\right) & {\rm \ if \ } v \in\overline{B}^4_2(0) \setminus  \overline{B}^4_1(0). \\
\end{array}
\right.
\end{eqnarray*}
 \subsection{Remark} \label{decreasing.mass} Notice that $${\bf M}(\tilde{C}(v,t)) \leq {\bf M}\left(C\left(\frac{v}{|v|},t\right)\right)$$ if $v\in \overline{B}^4_2(0) \setminus  \overline{B}^4_1(0)$ and $t\in [0,1]$.

 
\section{Min-max family}

In this section we construct the continuous map $\Phi:I^5 \rightarrow \mathcal{Z}_2(S^3)$ to which we apply Almgren-Pitts Min-Max Theory as described in \cite{marques-neves}.

Choose an orientation-preserving  homeomorphism $f:I^4 \rightarrow \overline{B}^4_2(0)$.

\subsection{Definition}\label{Fi.family} The {\bf min-max family} of  $(\gamma_1,\gamma_2)$ is the map $\Phi:I^5 \rightarrow \mathcal{Z}_2(S^3)$ given by
$$
\Phi(x,t) = \tilde{C}(f(x),t),
$$
where $x\in I^4$ and $t\in [0,1]$.

The properties of $\Phi$ that are important for our proof are collected in the next theorem.  If $T$ is an integral 2-current, we denote by $|T|$ the integral 2-varifold obtained from $T$ by forgetting orientations.  We denote by $\mathcal{T}$ the set of all  unoriented totally geodesic spheres, which is  homeomorphic to $\RP^3$. 

\subsection{Theorem}\label{modified.family}
{\em Let $(\gamma_1,\gamma_2)$ be a 2-component link in $S^3$ with ${\rm lk}(\gamma_1,\gamma_2)=-1$. The map $$\Phi:I^5 \rightarrow \mathcal{Z}_2(S^3)$$ satisfies  the following properties:
\begin{enumerate}
\item [(a)] $\Phi$ is continuous with respect to the flat topology of currents;
\item [(b)] $\Phi(I^4 \times \{0\})=\Phi(I^4 \times \{1\}) = \{0\};$
\item [(c)] $$
\sup\{{\bf M}(\Phi(x)):x\in I^5\}\leq E(\gamma_1,\gamma_2);
$$
\item[(d)] $\Phi(x,t)=\partial B_{r(t)}\left(\frac{f(x)}{|f(x)|}\right)$ for every $(x,t) \in \partial I^4 \times I$; 
\item[(e)] {the map $|\Phi|:\partial I^4\times\{1/2\}\rightarrow \mathcal{T}$ defined by
$$|\Phi|(x,1/2)= |\Phi(x,1/2)|=\left|\partial B_{\pi/2}\left(\frac{f(x)}{|f(x)|}\right) \right| \in \mathcal{T}$$
has
$$|\Phi|_{*}([\partial I^4\times\{1/2\}])=2\in H_3(\RP^3,\Z);$$ }
\item [(f)] $$
\lim_{r\to 0}\sup \left\{{\bf M}(\Phi(x)\llcorner B_r(q))\,:\, q\in S^3, x\in I^5\right\}=0. 
$$
\end{enumerate}
}

\begin{proof}
Property (a) follows from  Proposition \ref{equicontinuous.extension}, Proposition \ref{mapa.K} (i), Proposition \ref{mapa.K} (ii) and  the definition of $\tilde C$.

Property (b) follows from Proposition \ref{equicontinuous.extension} (i) and Remark \ref{decreasing.mass}.

Property (c) follows from Lemma \ref{canonical.area.bound}, Remark \ref{decreasing.mass}, and lower semicontinuity of mass under flat metric convergence. 

Property (d) follows from the fact that if $x\in \partial I^4$, then $|f(x)|=2$ and thus we have from Proposition \ref{mapa.K} (iii) that
$$\Phi(x,t)=  \tilde{C}(f(x),t)=C\left(\frac{f(x)}{|f(x)|}, t,1\right)= \partial B_{r(t)}\left(\frac{f(x)}{|f(x)|}\right).$$

Property (e) follows from the fact that $r(1/2)=\pi/2$, $f:\partial I^4\rightarrow \partial \overline B^4_2(0)$ is an orientation-preserving homeomorphism, and under the natural identification between $\mathcal{T}$ and $\RP^3$ we have
$$|\Phi|:\partial I^4\times\{1/2\}\rightarrow \RP^3,\quad |\Phi|(x,1/2)=\left\{ \frac{f(x)}{|f(x)|}, -\frac{f(x)}{|f(x)|}\right\}.$$

Property (f) follows from Proposition \ref{equicontinuous.extension} (iv) and Proposition \ref{mapa.K} (v).
\end{proof}


\section{Proof of the Main Theorem}

Since the energy $E$ is conformally invariant, we can assume the link $(\gamma_1,\gamma_2)$ is contained
in $S^3$. We can also assume ${\rm lk}(\gamma_1,\gamma_2) =-1$ by appropriately choosing the orientations 
of the curves.

Consider the min-max family $\Phi$ of $(\gamma_1,\gamma_2)$, as in Definition \ref{Fi.family}.  The idea is to apply the theory developed in \cite{marques-neves}  to the $5$-parameter family $\Phi$ and so we now argue that $\Phi$ satisfies all hypotheses necessary to apply  
Corollary 9.2 of \cite{marques-neves}.

 
Conditions {$(A_5)-(A_7)$} of \cite[Section 9]{marques-neves} are met due to  Theorem \ref{modified.family}.  Conditions $(A_0)-(A_3)$ of \cite[Section 8]{marques-neves} are also met due to Theorem \ref{modified.family}. 
The condition $(A_4)$ is there  to ensure that the continuous  path
$$\gamma:[0,1]\rightarrow\mathcal{Z}_2(S^3), \quad \gamma(t)=\Phi(1/2,1/2,1/2,1/2,t),$$
represents a non-trivial element of $\pi_1(\mathcal{Z}_2(S^3;\mathcal{F}),\{0\})$  (see \cite[Section 7]{marques-neves} for definition). Because $\Phi$ is continuous in the flat metric, we have that $\gamma$ and $t\mapsto\Phi(x,t)$, with $x\in \partial I^4$, represent the same element in $\pi_1(\mathcal{Z}_2(S^3;\mathcal{F}),\{0\})$. Note that $t\mapsto\Phi(x,t)$ is the standard sweepout of $S^3$ by round spheres, and thus it is non-trivial in $\pi_1(\mathcal{Z}_2(S^3;\mathcal{F}),\{0\})$ due to \cite[Theorem 7.1]{almgren}. Hence we checked that $\gamma$ defines a non-trivial element in $\pi_1(\mathcal{Z}_2(S^3;\mathcal{F}),\{0\})$ without using condition $(A_4)$. {Finally, we have 
$$\sup\{{\bf M}(\Phi(x)):x\in I^5\}\leq E(\gamma_1,\gamma_2)<8\pi$$
because we are assuming that ${\rm lk}(\gamma_1,\gamma_2) =-1$.} 

Therefore we can apply Corollary 9.2 of \cite{marques-neves}  to conclude the existence of a smooth embedded minimal surface $\Sigma \subset S^3$ with genus $g\geq 1$ and such that
$$
{\rm area}(\Sigma)\leq \sup\{{\bf M}(\Phi(x)):x\in I^5\}.
$$

Theorem B of \cite{marques-neves} gives that  ${\rm area}( \Sigma)\geq 2\pi^2$ and so Theorem \ref{modified.family} (c) implies that
\begin{equation}\label{energy.inequality}2\pi^2\leq {\rm area}(\Sigma) \leq \sup\{{\bf M}(\Phi(x)):x\in I^5\}\leq E(\gamma_1,\gamma_2).
\end{equation}

 We have proved that $E(\gamma_1,\gamma_2) \geq 2\pi^2$. 
 
 Suppose now $E(\gamma_1,\gamma_2)=2\pi^2$.  

\subsection{Lemma}\textit{$\sup\, \{{\bf M}(C(v,z)):(v,z) \in B^4 \times (0,1)\} = 2\pi^2$.}
\begin{proof}
It follows from (\ref{energy.inequality}) and Definition \ref{Fi.family} that 
\begin{eqnarray*}
&&\sup\, \{{\bf M}(\tilde{C}(v,z)):(v,z) \in \overline{B}^4_2(0) \times [0,1]\}  \\
&& \hspace{3cm}=\sup\{{\bf M}(\Phi(x)):x\in I^5\}=2\pi^2.
\end{eqnarray*}

But from Remark \ref{decreasing.mass}, we get that
\begin{eqnarray*}
&&\sup\, \{{\bf M}(\tilde{C}(v,z)):(v,z) \in \overline{B}^4_2(0) \times [0,1]\} \\
&&\hspace{3cm}= \sup\, \{{\bf M}(C(v,z)):(v,z) \in \overline{B}^4 \times [0,1]\}.
\end{eqnarray*}
Since the mass functional is lower semicontinuous, we also get
\begin{eqnarray*}
&&\sup\, \{{\bf M}(C(v,z)):(v,z) \in \overline{B}^4  \times [0,1]\} \\
&&\hspace{3cm}= \sup\, \{{\bf M}(C(v,z)):(v,z) \in B^4 \times (0,1)\}
\end{eqnarray*}
and the lemma follows.
\end{proof}

\subsection{Lemma}\label{existence.vz}\textit{There exists $(v,z) \in \overline{B}^4 \times [0,1]$ such that 
\begin{eqnarray*}
\int_{(S^1\times S^1) \setminus S_v} |{\rm Jac\, }g_{(v,z)}| \, ds\, dt=E(\gamma_1,\gamma_2)=2\pi^2.
\end{eqnarray*}
}
\begin{proof}
	Let $(v_i,z_i) \in B^4 \times (0,1)$ be a sequence such that ${\bf M}(C(v_i,z_i))\to2\pi^2$. 
We can assume that $(v_i,z_i)$ converges, as $i \rightarrow \infty$, to $(v,z) \in \overline{B}^4 \times [0,1]$. 

We know that
$$
{\bf M}(C(v_i,z_i)) \leq \int_{S^1\times S^1} |{\rm Jac\, }g_{(v_i,z_i)}| \, ds\, dt\leq E(\gamma_1,\gamma_2)=2\pi^2.
$$
Hence, by letting $i$ go to infinity, the lemma follows from Proposition \ref{equicontinuous.extension} (iii).
\end{proof}

\subsection{Lemma}\label{rigidity.lemma}\textit{Let $(v,z) \in \overline{B}^4 \times [0,1]$  be as in Lemma \ref{existence.vz}. Then $z=1/2$.
}
  \begin{proof}
 From Proposition \ref{equicontinuous.extension} (iii) we have $0<z<1$ (otherwise $|{\rm Jac\, }g_{(v,z)}|=0$).  From Lemma \ref{canonical.area.bound} we see that if $v\in B^4$ then
 \begin{multline*}
2\pi^2=\int_{S^1\times S^1}|{\rm Jac}\,g_{(v,z)}|ds\,dt \\
\leq   \int_{S^1\times S^1}\frac{a(v,z)|(F_{v}\circ \gamma_1)'(s)||(F_{v}\circ \gamma_2)'(t)|}{a(v,z)|F_{v}\circ \gamma_1(s)-F_{v}\circ \gamma_2(t)|^2+b(v,z)^2}ds\,dt\\
\leq E(\gamma_1,\gamma_2)=2\pi^2.
\end{multline*}
Thus $b(v,z)=0$ and  $z=1/2$.
 
 Suppose now $v\in S^3$. In this case we have
 \begin{align*}
 |F_v\circ\gamma_1-L_{(v,z)}\circ \gamma_2|^2 & =|(F_v\circ\gamma_1+v/2)-(F_v\circ \gamma_2+v/2)+b(z)v|^2\\
 & =|F_v\circ\gamma_1-F_v\circ \gamma_2|^2+b(z)^2.
 \end{align*}
 Thus from the expression for $g_{(v,z)}$ and Lemma \ref{jacobian.estimate} (i) we obtain that
 $$ |{\rm Jac\, }g_{(v,z)}|\leq \frac{|(F_v\circ\gamma_1)'||(F_v\circ\gamma_2)'|}{|F_v\circ\gamma_1-F_v\circ \gamma_2|^2+b(z)^2}.$$
From conformal invariance we have
 \begin{multline*}\int_{(S^1\times S^1)\setminus S_v}\frac{|(F_v\circ\gamma_1)'||(F_v\circ\gamma_2)'|}{|F_v\circ\gamma_1-F_v\circ \gamma_2|^2}\, ds\,dt=E(\gamma_1,\gamma_2)\\
 = \int_{(S^1\times S^1) \setminus S_v} |{\rm Jac\, }g_{(v,z)}| ds\, dt\\
 \leq \int_{(S^1\times S^1)\setminus S_v} \frac{|(F_v\circ\gamma_1)'||(F_v\circ\gamma_2)'|}{|F_v\circ\gamma_1-F_v\circ \gamma_2|^2+b(z)^2} ds\,dt.
 \end{multline*}
 This implies $b(z)=0$, i.e., $z=1/2$.
 \end{proof}
 
 We now prove the rigidity statement. Choose $(v,1/2)$ as  in the previous lemma. Then $E(\gamma_1,\gamma_2)$ is a minimizer of the M\"obius energy over all 2-component links with $lk(\gamma_1,\gamma_2)\neq 0$. By \cite[Section 7]{freedman-he-wang}, the minimizer must be smooth, since $\gamma_1$ is a geodesic with respect to a conformally invariant metric $ds_{E(\gamma_2)}^2$ on the complement of $\gamma_2(S^1)$, and vice-versa. 
 From Lemma \ref{jacobian.estimate} (i) we have that for all $(s,t)$
 $$\langle \hat{\gamma}_1'(s),\hat{\gamma}_2'(t)\rangle =\langle \hat{\gamma}_1'(s),\hat{\gamma}_1(s)-\hat{\gamma}_2(t)\rangle=\langle \hat{\gamma}_2'(t),\hat{\gamma}_1(s)-\hat{\gamma}_2(t)\rangle = 0,$$
where $\hat{\gamma}_i=F_v \circ \gamma_i$, $i=1,2$.  Therefore there exists $r_0>0$ such that $|\hat{\gamma}_1(s)-\hat{\gamma}_2(t)|^2=r_0^2$.
So $\hat{\gamma}_i(S^1)$ must be a round circle $i=1,2$. By \cite[Corollary 7.3]{freedman-he-wang}, $\hat{\gamma}_1(S^1)\cup\hat{\gamma}_2(S^1)$ must be
the standard Hopf link, up to M\"obius transformation. 

\bibliographystyle{amsbook}

\begin{thebibliography}{99}


\bibitem{almgren} 
F. Almgren, \textit{The homotopy groups of the integral cycle groups,} Topology  (1962), 257--299. 



\bibitem{burago-burago-ivanov}
D. Burago, Y. Burago, S. Ivanov,
\textit{A course in metric geometry,}
Graduate Studies in Mathematics 33, American Mathematical Socitey,  (2001).

\bibitem{freedman-he-wang}
M. Freedman, Z-X. He, Z. Wang, 
\textit{M\"{o}bius energy of knots and unknots,}
Ann. of Math. (2) 139 (1994), no. 1, 1--50. 

\bibitem{He02}
Zheng-Xu He, \emph{On the minimizers of the {M}\"obius cross energy of links},
  Experiment. Math. \textbf{11} (2002), no.~2, 244--248.


\bibitem{kim-kusner}
D. Kim and R. Kusner,
\textit{Torus knots extremizing the M\"{o}bius energy,}
Experiment. Math. Volume 2, Issue 1 (1993), 1-9. 

\bibitem{marques-neves} F. C. Marques and A. Neves,
\textit{Min-Max theory and the Willmore conjecture}, {Ann. of
Math. (2) 179 (2014), no. 2, 683--782.}

\bibitem{morgan} F. Morgan, {\it Geometric measure theory: a beginnerÕs guide, third edition},  Academic Press, 2000. 

\bibitem{o'hara}
J. O'Hara, \textit{Energy of a knot,} Topology 30 (1991), 241--247.

\bibitem{simon}
L. Simon, \textit{Lectures on geometric measure theory,}
Proceedings of the Centre for Mathematical Analysis, Australian National University,  Canberra, (1983). vii+272 pp. 

\bibitem{willmore}
T. J. Willmore, 
\textit{Note on embedded surfaces,} An. Sti. Univ.``Al. I. Cuza'' Iasi Sect. I a Mat. (N.S.) 11B (1965) 493--496. 


\end{thebibliography}

\end{document}